\documentclass[11pt]{article}%{elsart}
\usepackage{amsthm, amsmath, amssymb,epsfig}
\newcommand{\oa}{\overline{a}}
\newcommand{\ox}{\overline{x}}
\newcommand{\gb}{\operatorname{gb}}
\newcommand{\lpp}{\operatorname{lpp}}
\newcommand{\lc}{\operatorname{lc}}
\newcommand{\lcm}{\operatorname{lcm}}
\newcommand{\lex}{\operatorname{lex}}
\newcommand{\V}{\mathbb V}
\newcommand{\I}{\mathbb I}
\newcommand{\hs}{\hspace*{3mm}}

\newtheorem{thm}{Theorem}%[section]
\newtheorem{lem}[thm]{Lemma}
\newtheorem{prop}[thm]{Proposition}

\newtheorem{conj}[thm]{Conjecture}

\theoremstyle{definition}
\newtheorem{defn}[thm]{Definition}

\theoremstyle{remark}

\title{Improving DISPGB Algorithm Using the Discriminant Ideal
 \footnote{Work partially supported by the
    Ministerio de Ciencia y Tecnolog\'{\i}a under project
    BFM2003-00368, and by the Generalitat de
    Catalunya under project 2005 SGR 00692}}
\author{Montserrat Manubens, Antonio Montes\\
Departament de Matem\`atica Aplicada 2,\\ Universitat
Polit\`ecnica de Catalunya, Spain.\\
e-mail: antonio.montes@upc.edu\\
http://www-ma2.upc.edu/$\sim$montes}
\date{November, 2004}

\begin{document}
%\titlepage{
%\null \vspace{2.5truecm} \thispagestyle{empty}
%\thispagestyle{empty}
%\begin{center}
%\textbf{\large Improving DISPGB Algorithm \\ Using the Discriminant Ideal}\\[3mm]
%\textsc{Montserrat Manubens, Antonio Montes}
% \vspace{0.8truecm}
%
%MA2--IR--04--00015
%\end{center}
%}
%\newpage
%\null \thispagestyle{empty}
%\newpage
\setcounter{page}{1}

\maketitle

%\begin{frontmatter}

%\title{Improving DISPGB Algorithm Using the Discriminant Ideal}
%%\shorttitle{Improving DISPGB Using Discriminant Ideal}
%%\shortauthor{M.~Manubens, A.~Montes}%

%\author{Montserrat Manubens\thanksref{label1}},
%\author{Antonio Montes\thanksref{label1}}
%\corauth[cor]{Corresponding author: Antonio Montes.}
%
%\address{Departament de Matem\`{a}tica Aplicada II. \\Universitat Polit\`{e}cnica de Catalunya. Spain.}
%
%\thanks[label1]{Work partially supported by the
%    Ministerio de Ciencia y Tecnolog\'{\i}a under project
%    BFM2003-00368, and by the Generalitat de
%    Catalunya under project 2001 SGR 00224}
%
%\ead{montserrat.manubens@upc.edu}
%\ead{antonio.montes@upc.edu}
%\ead[url]{http://www-ma2.upc.edu/$\sim$montes}

\begin{abstract}
  In 1992, V. Weispfenning proved the existence of Comprehensive
  Gr\"obner Bases (CGB) and gave an algorithm to compute one. That
  algorithm was not very efficient and not canonical. Using his suggestions,
  A. Montes obtained in 2002 a more efficient algorithm
  (DISPGB) for Discussing Parametric Gr\"{o}bner Bases. Inspired in
  its philosophy,
  V. Weispfenning defined, in 2002, how to obtain a Canonical Comprehensive
  Gr\"obner Basis (CCGB) for parametric polynomial ideals, and provided a constructive method.

  In this paper we use Weispfenning's CCGB ideas to make
  substantial improvements on Montes DISPGB algorithm. It now includes rewriting
  of the discussion tree using the Discriminant Ideal and provides a
  compact and effective discussion.
  We also describe the new algorithms in the DPGB library containing
  the improved DISPGB as well as new routines to check whether a given basis is a CGB
  or not, and to obtain a CGB. Examples and tests are also
  provided.
\end{abstract}



\section{Introduction}
Let $R=k[\oa]$ be the polynomial ring in the parameters
$\oa=a_1,\dots,a_m$ over the field $k$, and $S=R[\ox]$ the
polynomial ring over $R$ in the set of variables
$\ox=x_1,\dots,x_n$. Let $\succ_{\ox}$ denote a monomial order wrt
the variables $\ox$, $\succ_{\oa}$ a monomial order wrt the
parameters $\oa$ and $\succ_{\ox\oa}$ the product order. The
problem we deal with consist of solving and discussing parametric
polynomial systems in $S$.

Since Gr\"obner bases were introduced various approaches have been
developed for this problem. The most relevant ones are:

\begin{itemize}
 \item Comprehensive Gr\"obner Bases (CGB)~\cite{We92}.
 \item Specific Linear Algebra Tools for Parametric Linear
 systems~\cite{Si92}.
 \item Dynamic Evaluation~\cite{Du95}.
 \item Newton Algorithm with Branch and Prune
 Approach~\cite{HeMcKa97}.
 \item Triangular Sets \cite{Mor97}.
 \item Specialization through Hilbert Functions~\cite{GoTrZa00}.
 \item DISPGB Algorithm~\cite{Mo02}.
 \item Alternative Comprehensive Gr\"obner Bases (ACGB)~\cite{SaSu03}.
 \item Canonical Comprehensive Gr\"obner Bases (CCGB)~\cite{We03}.
\end{itemize}

This paper describes some improvements made on {\tt DISPGB}.
Trying to solve some of the examples given in the references cited
above using the improved {\tt DISPGB} has been an interesting
challenge (see section \ref{Examples}).

In \cite{We92}, Professor Volker Weispfenning proved the existence
of a Comprehensive Gr\"obner Basis CGB wrt $\succ_{\ox}$ for any
ideal $I\subset S$ such that for every specialization of the
parameters $\sigma_{\oa}:R\rightarrow K'$ extended to
$R[\ox]\rightarrow K'[\ox]$, $\sigma_{\oa}(\hbox{CGB})$ is a
Gr\"obner basis of the specialized ideal $\sigma_{\oa}( I ) $. He
also provided an algorithm to compute it. There are two known
implementations of this algorithm \cite{Pe94,Sc91}.

In \cite{Mo95} and \cite{Mo98}, A. Montes used classical Gr\"obner
bases theory to study the load-flow problem in electrical
networks. V. Weispfenning recommended him to use the Comprehensive
Gr\"obner Basis algorithm~\cite{We92,Pe94} for this problem. The
use of CGB in the load-flow problem provided interesting
information over the parameters, but was rather complicated and
not very efficient. Moreover, it was not canonical, i.e. it was
algorithm depending.

In Montes~\cite{Mo02} provided a more efficient algorithm ({\tt
DISPGB}) to Discuss Parametric Gr\"obner Bases, but it was still
non-canonical. {\tt DISPGB} produces a set of non-faithful,
canonically reduced Gr\"{o}bner bases (Gr\"obner system) in a
dichotomic discussion tree whose branches depend on the
cancellation of some polynomials in $R$. The ideas in {\tt DISPGB}
however, inspired V. Weispfenning in \cite{We02, We03} to prove
the existence of a Canonical Comprehensive Gr\"obner Basis (CCGB)
as well as to give a method to obtain one.

The main idea for building up the canonical tree is the obtention
of an ideal $J\subset R$, structurally associated to the ideal
$I\subset S$ and the order $\succ_{\ox}$, which clearly separates
the essential specializations not included in the generic case.
Let us denote $J$ as the {\em Weispfenning's discriminant ideal}
of $(I,\succ_{\ox})$. In the new Weispfenning's algorithm, $J$
must be computed at the beginning of the discussion using a
relatively time-consuming method. The discriminant ideal was one
of the lacks of the old {\tt DISPGB} and an insufficient
alternative algorithm {\tt GENCASE} was provided.

In this paper we obtain, following Weispfenning, a {\em
discriminant ideal} denoted as $N$, which can be determined from
the data obtained after building the {\tt DISPGB} tree using a
less time consuming algorithm and, moreover, we prove that
$J\subset N$. We conjecture that $J=N$. We have verified it in
more than twenty different examples, and no counter-example has
been found. The ideal $N$ allows to rewrite the tree getting a
strictly better discussion.

We also prove that for a large set of parametric polynomial ideals
(at least for all prime ideals $I$) the discriminant ideal is
principal and in this case we have a unique {\em discriminant
polynomial} to distinguish the generic case from the essential
specializations. All the theoretical results commented above are
detailed in section \ref{DiscrimIdeal}.

In section \ref{ImprovedDISPGB}, we describe the improvements
introduced in the algorithms. We have made a complete revision to
the old release simplifying the algorithm and highly increasing
its speed. New routines {\tt CANSPEC} and {\tt PNORMALFORM} which
perform semi-canonical specifications of specializations and
reductions of polynomials are given. The algorithm has been
completely rewritten and the flow control has been simplified.
Further reductions of the tree, eliminating similar brother
terminal vertices, have been performed using algorithm {\tt
COMPACTVERT}.

Following P. Gianni~\cite{Gi87}, we are interested in guessing
whether some basis of $I$ is a comprehensive Gr\"obner basis or
not, in particular for the reduced Gr\"obner basis of $I$ wrt the
product order $\succ_{\ox\oa}$. We give, in section \ref{CGBalg},
a simple algorithm {\tt ISCGB} which uses the {\tt DISPGB} output
tree to answer that question. We also give an algorithm {\tt
PREIMAGE} to compute a faithful pre-image of the non-faithful
specialized polynomials from the reduced bases. This allows to
construct a CGB. It will be interesting to compare our CGB with
Weispfenning's CCGB when implemented.

Finally, in section \ref{Examples}, we give two illustrative
examples and a table of benchmarks for {\tt DISPGB} applied to
several parametric systems from which the power of the algorithm
is clearly shown.

It is stated in the same section that the new {\tt
DISPGB}\footnote{Release 2.3 of the library {\tt DPGB}, actually
implemented in {\em Maple} and available at the site
http://www-ma2.upc.edu/$\sim$montes/} algorithm is efficient and
provides a compact discussion of parametric systems of polynomial
equations. An incipient version of it was presented in
\cite{MaMo04}.

\section{Generic Case, Discriminant Ideal and Special
Cases}\label{DiscrimIdeal}

Let $K=k(\oa)$ be the quotient field of $R$ and $IK$ the ideal $I$
extended to the coefficient field $K$. Consider
$G=\gb(IK,\succ_{\ox})$, the reduced Gr\"obner basis of $I  K$ wrt
$\succ_{\ox}$. As $K$ is a field, $G$ can be computed through the
ordinary Buchberger algorithm. The polynomials in $G$ have leading
coefficient 1. With this normalization $g$ can have denominators
in $R$. Let $d_g\in R$ be the least common multiple of the
denominators of $g$. To obtain a polynomial in $S$ corresponding
to $g$ it suffices to multiply $g$ by $d_g$. Following
Weispfenning~\cite{We02,We03}, for each $g\in G$ we can obtain a
{\em minimal lifting} of $g$, $a_gg$, such that $a_gg\in I$ and
$a_g\in R$ is minimal wrt $\succ_{\oa}$. Doing this for all $g\in
G$ we obtain $G'$, a minimal lifting of $G$ which Weispfenning
calls the {\em generic Gr\"obner basis} of $(I,\succ_{\ox})$. Of
course, $d_g \mid a_g$. We will use a sub-lifting of $G$, $G''=\{
d_gg \ : \ g \in G \} \subset S$, and this will be our {\em
generic case basis} because it is simpler to compute and
corresponds to our standard form of reducing polynomials, as it
will be seen in section \ref{ImprovedDISPGB}.

We call {\em singular specialization} a specialization $\sigma$
for which the set of $\lpp$ (leading power products) of the
reduced Gr\"obner basis of $\sigma(I)$ is not equal to the set of
$\lpp(G,\succ{\ox})$.

{\tt DISPGB} builds up a binary dichotomic tree
$T(I,\succ_{\ox},\succ_{\oa})$ branching at the vertices whenever
a decision about the cancellation of some $p\in R$ has been taken.
Each vertex $v \in T$ contains the pair $(G_v,\Sigma_v)$.
$\Sigma_v=(N_v,W_v)$ is the semi-canonical specification of the
specializations in $v$, where $N_v$ is the radical ideal of the
current assumed null conditions (from which all factors of
polynomials in $W_v$ have been dropped), and $W_v$ is the set of
irreducible polynomials (conveniently normalized and reduced by
$N_v$) of the current assumed non-null conditions. Considering
$W^{*}_v$ the multiplicatively closed set generated by $W_v$, then
$G_v\subset (W^{*}_v)^{-1}\left( K[\ox]/N_v \right)$ is the
reduced form of the basis of $\sigma(I)$ for the specification of
the specializations $\sigma \in \Sigma_v$. At a terminal vertex,
the basis $G_v$ is the reduced Gr\"obner basis of $\sigma(I)$, up
to normalization, for all specializations $\sigma \in \Sigma_v$.

 Weispfenning \cite{We02} introduces the
following ideal associated to each $g\in G$:
\[J_g= \{a\in R \ : \ ag \in I \} = d_g\,  ( I:d_g g) \bigcap R\]
the second formula being computable via ordinary Gr\"obner bases
techniques. Then the radical of their intersection
$J=\sqrt{\bigcap_{g\in G} J_g}$ is used to distinguish the generic
case in the algorithm. We call $J$ the {\em Weispfenning's
discriminant ideal}. A specialization $\sigma$ is said to be {\em
essential} (for $I,\succ_{\ox}$) if $J_g \subseteq \ker(\sigma)$
for some $g\in G$.

V. Weispfenning proves the following two theorems:

\hs\hs \begin{minipage}[t]{4.8in}
\begin{itemize}
  \item [W1:] $J=\bigcap \ \{ \ker(\sigma) \ : \ \sigma \hbox{ is essential }
\}.$
  \item [W2:]  Let $\sigma$ be an inessential specialization.
  Then
 \begin{itemize}
  \item [(i)] $\sigma(G)$ is defined for every $g\in G$ and
  $\lpp(\sigma(g),\succ_{\ox})=\lpp(g,\succ_{\ox})$.
  \item [(ii)]$\sigma(G)$ is the reduced Gr\"obner basis of the ideal
  $\sigma(I)$.
 \end{itemize}
\end{itemize}
\end{minipage}

In the {\tt DISPGB} tree $T(I,\succ_{\ox},\succ_{\oa})$
specializations are grouped into disjoint final cases $i$ by the
specification $\Sigma_i$, and for all specializations in
$\Sigma_i$ the reduced Gr\"obner bases have the same set of $\lpp$
wrt $\succ_{\ox}$.

Let $1 \le i \le k$ number the terminal vertices. We call {\em
singular cases} the final cases for which $\lpp(G_i,\succ_{\ox})
\ne \lpp(G,\succ_{\ox})$. Let $A$ be the set of indexes of the
singular cases:
\[A=\{   1 \le i \le k \ : \ \lpp(G_i,\succ_{\ox}) \ne
\lpp(G,\succ_{\ox}) \}. \] We denote $\V(I)$ the variety of $I$
and $\I(V)$ the ideal of the variety $V$.
 The tree, being dichotomic, provides a partition of $(K')^m$ into disjoint
sets of specifications, and thus
\[(K')^m=\bigcup_{i=1}^k \left( \V(N_i) \setminus  \bigcup_{w\in W_i} \V(w)  \right)= U_s \ \bigcup \ U_g, \]
where $U_s$ is the set of points $\oa \in (K')^m$ corresponding to
singular specifications, i.e.
\[U_s(I,\succ_{\ox})=\{\oa \in (K')^m \ : \ \sigma_{\oa} \hbox{ is singular } \} =
\bigcup_{i \in A} \left( \V(N_i) \setminus  \bigcup_{w\in W_i}
\V(w)  \right).
\]
\begin{thm} \label{NinterNi} Let us call $N(I,\succ_{\ox})=\I(U_s)$ the {\em discriminant
ideal}. Then $$N(I,\succ_{\ox})=\bigcap_{i\in A} N_i.$$
\end{thm}
This theorem allows to compute $N$ from the output of {\tt
BUILDTREE}, i.e. the first tree construction in {\tt DISPGB}. (See
section \ref{ImprovedDISPGB}).
\begin{proof} We prove both inclusions:
\begin{itemize}
\item[$\subseteq$:] $f(\oa)=0$ for all $f\in N=\I(U_s)$ and $\oa \in U_s$.
Thus $\sigma_{\oa}(f)=0$ for all $\oa \in U_s$. Taking now $\oa$
such that $\sigma_{\oa} \in \Sigma_i$ this implies that $f \in
N_i$. As this can be done for all $i \in A$, it follows that $N
\subseteq \bigcap_{i \in A} N_i$.
\item[$\supseteq$:] For all $f\in \bigcap_{i\in A} N_i$ and all
$\oa\in U_s$ there exists $i\in A$ such that $\sigma_{\oa} \in
\Sigma_i$ and, of course, $f\in N_i$. Thus $\sigma_{\oa}(f)=0$,
i.e. $f(\oa)=0$ for all $\oa \in U_s$. Thus $f \in \I(U_s)=N$.
\end{itemize}
\end{proof}
Before proving the next theorem we need the following

\begin{lem} \label{singisessen} Any singular specialization is
essential.
\end{lem}
\begin{proof} Let $\sigma_{\oa}$ be a singular specialization. If it were not
essential, by Weispfenning theorem (W2), then the reduced
Gr\"obner basis of $\sigma(I)$ would be the generic basis $G$, and
this contradicts the definition of singular specialization. Thus
$\sigma_{\oa}$ must be essential.
\end{proof}

\begin{thm} $J \subseteq N$.
\end{thm}
\begin{proof}  By Weispfenning's theorem (W1), if $f\in J$
then $f\in \ker(\sigma_{\oa})$ for all essential $\sigma_{\oa}$,
and thus $f(\oa)=0$. So, by lemma \ref{singisessen}, $f(\oa)=0$
for all singular $\sigma_{\oa}$. This implies that $f(\oa)=0$ for
all $i\in A$ and $\sigma_{\oa} \in \Sigma_i$ and thus $f\in
\sqrt{N_i}=N_i$. Finally, by proposition \ref{NinterNi}, $f\in N$.
\end{proof}

\begin{conj}\label{conject} We formulate two forms

\hs\hs\begin{minipage}[t]{4.8in}
\begin{itemize}
  \item [{\rm(}i{\rm)}] {\rm(}Strong conjecture{\rm)}. All essential specializations
  are singular.
  \item [{\rm(}ii{\rm)}] {\rm(}Weak conjecture{\rm)}. $J \supseteq N$.
\end{itemize}
\end{minipage}
\end{conj}
\begin{prop}
The strong formulation of conjecture \ref{conject} implies the
weak formulation.
\end{prop}
\begin{proof} If $f\in N$ then, for all $i\in A$, $f\in N_i$. Thus, $f(\oa)=0$
for all singular specialization $\sigma_{\oa}$ and, if the strong
form of the conjecture is true, then $f(\oa)=0$ also for all
$\sigma_{\oa}$ essential and thus $f\in \ker(\sigma_{\oa})$. So,
by Weispfenning's theorem (W1), $f\in J$.
\end{proof}
In any case, by definition $N$ is discriminant, i.e. for any
$\oa\not\in \V(N)$ the Gr\"obner basis of $\sigma_{\oa}(I)$ is
generic, and every singular specification is in $\V(N)$. Thus,
what we called minimal singular variety in \cite{Mo02} is
described by $\V(N)$. If the strong formulation of the conjecture
is true then every specialization $\sigma$, for which $N \subset
\ker(\sigma)$, is not only essential but also singular and thus
the corresponding set of $\lpp$ of its reduced Gr\"obner basis
cannot be generic.

We have tested our conjecture in more than twenty examples and we
have not found any counter-example of any of the two formulations.
Nevertheless the weak formulation is the most interesting one and
a failure of the strong formulation would not necessarily
invalidate the weak formulation.

In most cases Weispfenning's discriminant ideal $J$ is principal,
as states the following

\begin{thm} \label{principalsg}
If $I\subset S$ is a prime ideal and the generic Gr\"obner basis
$G$ wrt $\succ_{\ox}$ is not $[1]$, then the discriminant ideal
$J(I,\succ_{\ox})$ is principal and is generated by the radical of
the $\lcm$ of all the denominators of the polynomials in $G$.
\end{thm}
\begin{proof}
Take $g\in G$. We have $J_g= d_g\,  (  I  : d_gg ) \bigcap R$. If
$h\in J_g$ then $d_g \mid h$, as $d_gg$ has no common factor with
$d_g$. Thus $d_gg\,(h/d_g) \in I$. By hypothesis, $d_gg\ne 1$ and
$I$ is prime. So, as $h/d_g \in R$, we have $h/d_g \not\in I$.
Thus, necessarily $d_gg \in I$ and $d_g\in J_g$. As $d_g \mid h$
for all $h \in J_g$, it follows that $J_g=\langle d_g \rangle$ is
principal. As $J=\sqrt{\bigcap_{g \in G} J_g}$ is the intersection
of principal ideals, the proposition follows.
\end{proof}
Not only prime ideals have principal discriminant ideals as the
next example shows: Take \[I=\langle a x+y+z+b,x-1+a y+z+b,x+y+a
z+b\rangle.\] Computing the Gr\"obner basis of $I$ wrt
$\lex(x,y,z,a,b)$ one can see that
\[I=\langle (a+2) z+b, y-z, x+y+a z+b \rangle \cap \langle
a-1,x+y+a z+b \rangle \] and $I$ is not prime. The generic
Gr\"obner basis wrt $\lex(x,y,z)$ is, in this case,
$G=[z+b/(a+2),y+b/(a+2),x+b/(a+2)]$. Thus $d_g=a+2$ for each $g\in
G$. For this example it is easy to compute $J=\langle (a+2)(a-1)
\rangle$ which is still principal even if $I$ is not prime and has
a prime component with generic Gr\"obner basis $[1]$.

It would be interesting to characterize which ideals $I\subset S$
have principal discriminant and which do not. But it is now clear
that in the most interesting cases we have principal
discriminants. This gives a new insight into our concept of
singular variety used in the algorithm~\cite{Mo02} in order to
understand the parallelism and differences between the new
Weispfenning's algorithm \cite{We02,We03} and {\tt DISPGB}, and
allows us to improve the old algorithm.

Under that perspective, we have completely revised ~\cite{Mo02}
and obtained a much more efficient and compact discussion. An
intermediate version was presented in~\cite{MaMo04}. We shall
describe now the improvements introduced in the new {\tt DPGB}
library and refer to \cite{Mo02}, where the old {\tt DPGB} is
described, for all unexplained details.

\begin{table}
\begin{minipage}[t]{5.2in}
\begin{tabular}{||p{1in}|p{1.4in}|p{1.6in}|p{0.7in}||}
\hline\hline Routines of the old al\-go\-rithm & Routines of the new algorithm & Improvements & Obsolete routines \\
\hline\hline
 {\tt DISPGB}\newline \hs {\tt BRANCH} %(r)%\footnote{Recursive algorithm are indicated with (r)}
 \newline
 &  {\tt DISPGB} \newline \hspace*{1mm}{\tt BUILDTREE}  \newline \hspace*{1mm}{\tt
 DISCRIMINANTIDEAL}\newline \hspace*{1mm}{\tt REBUILDTREE}  \newline \hspace*{1mm}{\tt COMPACTVERT}
 \newline
 & {\tt BUILDTREE} replaces old {\tt BRANCH}.

 Current {\tt DISPGB} includes also re\-buil\-ding of the
 tree ({\tt REBUILDTREE}) and {\tt
 COMPACTVERT}.
 & {\tt GENCASE} \\
 \hline
 {\tt BRANCH}  \newline \hs {\tt NEWVERTEX} & {\tt BUILDTREE} & Better flow control,
 no incom\-pa\-ti\-ble branching. & {\tt BRANCH} \\
 \hline
 {\tt NEWCOND} & {\tt CONDTOBRANCH} & More robust, ensures no
 incompatible branches. & {\tt NEWCOND} \\
 \hline
 {\tt CANSPEC} & {\tt CANSPEC} & Uses radical ideal. More
 robust. & \\
 \hline
 \hs \hs - & {\tt PNORMALFORM} & Standard polynomial reduction wrt $\Sigma$. & \\
 \hline {\tt CONDPGB} & {\tt CONDPGB} & Uses {\tt
 CONDTOBRANCH} and Weis\-pfen\-ning's standard pair selection. &
 \\
 \hline
 \hs \hs - & {\tt DISCRIMINANTIDEAL} & Determines the discriminant ideal $N$. &\\
 \hline
 \hs \hs - & {\tt REBUILDTREE} & Rebuilds the tree starting the
 discussion with $N$.  & {\tt GENCASE}  \newline (external) \\
 \hline
 \hs \hs - & {\tt COMPACTVERT} & Drops brother terminal vertices with
 same $\lpp$ sets. & \\
 \hline\hline
\end{tabular}
\end{minipage}
\newline

\caption{\label{impdispgb}}
\end{table}
\section{Improved DISPGB Algorithm}\label{ImprovedDISPGB}
In this section we describe the improvements introduced in {\tt
DISPGB} algorithm. Table \ref{impdispgb} summarizes the basic
differences between old \cite{Mo02} and the new algorithms used in
it.

First, we have improved the construction of the discussion tree
$T(I,\succ_{\ox},\succ_{\oa})$ in order to have a simpler flow
control and to make it faster by avoiding un\-ne\-ces\-sa\-ry and
useless time-consuming computations. In the old algorithm this was
done by the recursive routine {\tt BRANCH} which was the unique
action of {\tt DISPGB}, but now it is done by {\tt BUILDTREE}. As
we explain later, it has been strongly reformed.

Then, {\tt DISCRIMINANTIDEAL} computes the discriminant ideal
$N=\bigcap_{i\in A} N_i$ which, as shown in section
\ref{DiscrimIdeal}, can be determined from {\tt BUILDTREE} output.

After that, {\tt DISPGB} calls {\tt REBUILDTREE}. This algorithm
builds a new tree setting the discriminant ideal $N$ at the top
vertex and the generic case at the first non-null vertex labelled
as $[1]$ (see figure \ref{robotarm} in section \ref{robot}). The
old tree is rebuilt under the first null vertex recomputing the
specifications and eliminating incompatible branches. The result
is a drastic reduction of branches in the new tree. In the old
{\tt DPGB} library, this work was partially done by the external
algorithm {\tt GENCASE} which has become useless.

To further compact the tree, a new algorithm {\tt COMPACTVERT} is
used. It summarizes brother terminal vertices with the same set of
$\lpp$ into their father vertex. {\tt COMPACTVERT} is called
before and after {\tt REBUILDTREE}. {\tt DISPGB} algorithm is
sketched in Table \ref{dispgbalg}.

\begin{table}
\fbox{\parbox[t]{5.2in}{
 $T \leftarrow \hbox{{\bf DISPGB}}(B,\succ_{\ox}, \succ_{\oa})$ \newline
 {\tt Input}: \newline
 \hs \hs $B \subseteq R[\oa][\ox]$\ :  basis of $I$, \newline
 \hs \hs $\succ_{\ox}$, $\succ_{\oa}$\ :  termorders wrt the
 variables $\ox$ and the parameters $\oa$
 respectively. \newline
 {\tt Output}: \newline
 \hs \hs $T$: table with binary tree structure, containing $(G_v,\Sigma_v)$ at vertex $v$ \newline
 BEGIN \newline
 \hs $T:=\phi$,   \ \# global variable\newline
 \hs $v:=[\ ]$ \ \# (label of the current vertex) \newline
 \hs $\Sigma:=([\ ], \phi)$ \ \# (current specification) \newline
 \hs $\hbox{{\tt BUILDTREE}}(v,B,\Sigma)$ \# (recursive, stores the computations in $T$)\newline
 \hs $N:=\hbox{{\tt DISCRIMINANTIDEAL}}(T)$ \newline
 \hs $\hbox{{\tt COMPACTVERT}}(T)$ \# (compacts $T$) \newline
 \hs $\hbox{{\tt REBUILDTREE}}(T,N)$ \# (rebuilds $T$) \newline
 \hs $\hbox{{\tt COMPACTVERT}}(T)$ \# (compacts $T$) \newline
 END
}} %end parbox, end fbox
\newline

\caption{\label{dispgbalg}}
\end{table}

\subsection{Building up the Discussion Tree: {\tt
BUILDTREE}.}\label{Buildtree}

We have simplified the flow control from the ancient {\tt DISPGB}
and dropped useless operations. Now all the hard work of the
discussion is done by the recursive algorithm {\tt BUILDTREE}
which replaces the old {\tt BRANCH} routine and makes {\tt
NEWVERTEX} useless. The obtained discussion is equivalent to the
one given by the old {\tt DISPGB}, but now is more compact.

It computes the discussion tree faster than the old one because
now it assembles the discussion over the coefficients of the
current basis in one single algorithm, avoiding unnecessary
branching and useless computations.

Given $B$, a set of polynomials generating the current ideal, {\tt
BUILDTREE} takes the current basis $B_v$ at vertex $v$,
specialized wrt the current reduced specification $\Sigma_v =
(N_v,W_v)$, builds a binary tree $T$ containing the discussion
under vertex $v$, and stores all the data at the vertices of $T$.
It is a recursive algorithm and substitutes the old {\tt BRANCH}
and {\tt NEWVERTEX}. See table \ref{buildtreealg}.

Theorem 16 in \cite{Mo02} still applies to the reformed {\tt
BUILDTREE}, thus we can assert the correctness and finiteness of
the algorithm.

\begin{table}
\fbox{
\begin{minipage}[t]{5.2in}
 \noindent{\bf BUILDTREE}$(v,B,\Sigma)$ \newline
 {\tt Input}: \newline
 \hs\hs $v$, the label of the current vertex, \newline
 \hs\hs $B \subseteq R[\bar{a}][\bar{x}]$, the current basis, \newline
 \hs\hs $\Sigma=(N,W)$ the current reduced specification. \newline
 {\tt Output}:
 No output, but the data are stored in the global tree variable
 $T$.\newline
 BEGIN \newline
 \hs $c_f:= \hbox{false}$  \newline
 \hs ($c_b,c_d,G,\Sigma_0,\Sigma_1$):={\tt CONDTOBRANCH}($B,\Sigma$) \newline
 \hs IF $c_d$ THEN \ \# ($c_d$ is true if all $\lc(g)$, $g\in G$ are decided non-null, false otherwise) \newline
 \hs\hs ($c_b,c_f,G,\Sigma_0,\Sigma_1$):={\tt CONDPGB}($G,\Sigma$) \newline
 \hs END IF \newline
 \hs $T_v:=(G,\Sigma)$ \ \# (Store data in the global tree variable $T$) \newline
 \hs IF $c_f$ THEN \ \# ($c_f$ is true if the new vertex is terminal, false otherwise)\newline
 \hs\hs RETURN() \newline
 \hs ELSE \newline
 \hs\hs IF $c_b$ THEN \ \# ($c_b$ is true if null and non-null conditions are both compatibles) \newline
 \hs\hs\hs {\tt BUILDTREE}$((v,0),G,\Sigma_0)$ \newline
 \hs\hs\hs {\tt BUILDTREE}$((v,1),G,\Sigma_1)$ \newline
 \hs\hs ELSE \newline
 \hs\hs\hs {\tt BUILDTREE}$(v,G,\Sigma_1)$ \ \# (and {\tt BUILDTREE} continues
  in the same vertex)~\footnote{In this case, if {\tt CONDPGB} has already started then the list of
  known $S$-polynomials reducing to 0 can be kept.}  \newline
 \hs\hs END IF \newline
 \hs END IF \newline
 END \newline
\end{minipage}
}%end fbox
\newline
\caption{\label{buildtreealg}}
\end{table}

The most important algorithms used by {\tt BUILDTREE} are
commented below.

The algorithm {\tt CONDTOBRANCH} substitutes the old {\tt
NEWCOND}. It is used each time that {\tt BUILDTREE} is recursively
called and also inside {\tt CONDPGB}, applying it to each new
not-reducing-to-zero $S$-polynomial. This prevents Buchberger
algorithm from stopping and saves incompatible branches.

Each time we need to know whether a given polynomial $f\in R$
--\,for example the $\lc$ (leading coefficient) of a new
$S$-polynomial\,-- is zero or not for a given specification, we
will reduce it by $\Sigma=(N,W)$ using {\tt PNORMALFORM} and then
test whether the remainder is compatible or not with taking it
null and non-null for each of the specifications using {\tt
CANSPEC}. The whole task is done by {\tt CONDTOBRANCH}. See table
\ref{condtobranchalg}.

\begin{table}
\fbox{\parbox[t]{5.28in}{
 \noindent$(c_b,c_d,G,\Sigma_0,\Sigma_1)\leftarrow \hbox{{\bf CONDTOBRANCH}}(B,\Sigma)$ \newline
 {\tt Input}: \newline
 \hs\hs $B \subseteq R[\bar{a}][\bar{x}]$, the current basis \newline
 \hs\hs $\Sigma=(N,W)$ a reduced specification. \newline
 {\tt Output}: \
   \newline
 \hs\hs $G$ is $B$ reduced wrt $\Sigma$, \newline
 \hs\hs $\Sigma_1$ is the reduced specification for the not null branch\newline
 \hs\hs $\Sigma_0$ is the reduced specification for the null branch\newline
 \hs\hs $c_b$ is true whenever $\Sigma_0$ exists, and false otherwise. \newline
 \hs\hs $c_d$ is true if all $g\in G$ have $\lc(g)$ decided to not null,
 and false otherwise.\newline
 BEGIN \newline
 \hs\hs $G:=\hbox{{\tt PNORMALFORM}}(B, \Sigma)$ \newline
 \hs\hs IF there is $g \in G$ with $l_g=\lc(g)$ not yet decided to not
 null wrt $\Sigma$ THEN \newline
 \hs\hs\hs $c_d:=\hbox{false}$ \newline
 \hs\hs\hs $(t,\Sigma_1):=\hbox{{\tt CANSPEC}}( N_{\Sigma},
 W_{\Sigma} \bigcup \{l_g\})$ \newline
 \hs\hs\hs $(t,\Sigma_0):=\hbox{{\tt CANSPEC}}(\langle N_{\Sigma},l_g \rangle ,
 W_{\Sigma})$ \newline
 \hs\hs\hs IF $t$ THEN  $c_b:=\hbox{true}$ ELSE $c_b:=\hbox{false}$ ENDIF \newline
 \hs\hs ELSE \newline
 \hs\hs\hs $c_d:=\hbox{true}$ \newline
 \hs\hs ENDIF \newline
 END\newline
}}%end parbox, end fbox
\newline

\caption{\label{condtobranchalg}}
\end{table}

{\tt BUILDTREE} uses a Buchberger-like algorithm --\,{\tt CONDPGB}
(Conditional Parametric Gr\"obner Basis)\,-- taking the
specification into account and intending to determine a
specializing Gr\"obner basis. The basic improvements on {\tt
CONDPGB} in the new version are: the call to {\tt CONDTOBRANCH}
instead of old {\tt NEWCOND} and improving Buchberger algorithm by
considering Weispfenning's normal strategy of pair selection
\cite{BeWe93}. We do not detail these improvements.

{\tt CANSPEC} has also been modified.

At each vertex $v$ of the tree  a pair ($G_v,\Sigma_v$) is stored,
where $\Sigma_v=(N_v,W_v)$ is a specification of specializations.
This means that for all $\sigma \in \Sigma_v$ one has
$\sigma(N_v)=0$ and $\sigma(w)\not=0 ~\forall w\in W_v$. From the
geometric point of view, a given $\Sigma=(N,W)$ describes the set
of points $\V(N) \setminus ( \bigcup_{w\in W} \V(w)) \subseteq
(K')^m$. %of the affine space

By proposition 5 in \cite{Mo02}, one can see that $\Sigma=(N,W)$
and $\Sigma'=(\sqrt{N},W)$ describe equivalent specialization
sets. And, by proposition 7, the same happens with
$\widetilde{\Sigma}=(\widetilde{N},\widetilde{W})$, where
$\widetilde{N}$ has no factor laying in $W$ and is radical, and
$\widetilde{W}$ is the set of the irreducible factors of $W$ with
multiplicity one reduced modulus $\widetilde{N}$. So we choose the
following representative for the specifications describing
equivalent specialization sets:

\begin{defn}\label{canspecific}
  We call $\Sigma=(N,W)$ a {\em reduced specification of
  specializations} if it is a specification such that

  \hs\hs \begin{minipage}[t]{5.2in}  \begin{itemize}
   \item[{\rm(}i{\rm)}] $\langle N \rangle$ is a radical ideal, and $N=\gb(\langle N
   \rangle,\succ_{\oa})$,
   \item[{\rm(}ii{\rm)}] there is no factor of any polynomials in $\langle N \rangle$ laying within $W$,
   \item[{\rm(}iii{\rm)}] $W$ is a set of distinct irreducible polynomials not laying within $\langle N \rangle$,
   \item[{\rm(}iv{\rm)}] $\overline{W}^N = W$.
  \end{itemize}
  \end{minipage}
\end{defn}
We must note that the set $W$ is not uniquely determined, as there
exist infinitely many polynomials which cannot be null for a given
specification. For example, suppose that the current reduced
specification is $W=\{ a \},N=[a^2-1]$. The condition $a\ne 0$ is
compatible with $N$ but is redundant in this case. We can also add
to $W$ other polynomials like $a-2$. Thus there is no unique
reduced specification, but our choice is convenient enough. The
task of obtaining reduced specifications and testing compatibility
of the current null and non-null conditions is done by the
reformed {\tt CANSPEC}. See table \ref{canspecalg}.
\begin{table}
\fbox{\parbox[t]{5.2in}{

 \noindent $(t,\widetilde{\Sigma}) \leftarrow \hbox{{\bf CANSPEC}}(\Sigma)$ \newline
 {\tt Input}: $\Sigma=(N,W)$ a not necessarily
 reduced specification. \newline
 {\tt Output}: \newline
 \hs $t$: a boolean valued variable. \newline
 \hs $\widetilde{\Sigma}$: a reduced
 specification if $t=\hbox{true}$, and $\phi$ otherwise (in this case \newline
 \hs\hs  incompatible conditions have been found). \newline
 BEGIN \newline
 \hs $N_a:=N, ~N_b:=\sqrt{N}$ \newline
 \hs $W_a:=W, ~W_b:=$ the irreducible factors of $W$ without multiplicity and reduced wrt $N_a$; \newline
 \hs IF $\prod_{q\in W_b}{q} =0$ THEN RETURN(false,$\phi$) ENDIF \newline
 \hs WHILE ($N_a \not= N_b$ AND $W_a \not=W_b$) DO \newline
 \hs\hs $N_a:=\phi$ \newline
 \hs\hs FOR $p \in N_b$ DO  \newline
 \hs\hs\hs $p:=$ drop from $p$ all irreducible factors laying in $W_b$ \newline
 \hs\hs\hs IF $p=1$ THEN RETURN(false,$\phi$) ENDIF \newline
 \hs\hs\hs Add $p$ into $N_a$ \newline
 \hs\hs END FOR \newline
 \hs\hs $W_a:= W_b$ \newline
 \hs\hs $N_b:=\sqrt{N_a}$ \newline
 \hs\hs $W_b:=$ $\{$irreducible factors of $W_a$ without multiplicity and reduced wrt $N_b\}$ \newline
 \hs\hs IF $\prod_{q\in W_b}{q} =0$ THEN RETURN(false,$\phi$) ENDIF \newline
 \hs END WHILE \newline
 \hs $\widetilde{\Sigma}:=(N_a,W_a)$ \newline
 \hs RETURN(true, $\widetilde{\Sigma}$) \newline
 END \newline
}} %end parbox, end fbox
\newline

\caption{\label{canspecalg}}
\end{table}

\begin{prop}
Given any specification of specializations $\Sigma=(N,W)$, if {\tt
CANSPEC}\,{\rm(}$\Sigma${\rm)} returns $(t,\widetilde{\Sigma})$
with $t=\hbox{true}$, then $\widetilde{\Sigma}$ is a reduced
specification of $\Sigma$ computed in finitely many steps.
Otherwise it returns $t=\hbox{false}$ and $(N,W)$ are not
compatible conditions.
\end{prop}
 \begin{proof}
At the end of each step $N_a$ is a radical ideal, $W_a$ is a set
of irreducible polynomials with multiplicity one reduced wrt
$N_a$, so $\overline{W_a}^{N_a}=W_a$. So, $N_b$ is still radical
when the algorithm stops, as $N_b$ is built by dropping from $N_a$
all those factors laying in $W_a$. If the algorithm returns {\tt
true}, as at each completed step $(N_b,W_b)$ satisfies the
conditions of definition \ref{canspecific}, then the conditions
are compatible and $\widetilde{\Sigma}$ is a reduced specification
of specializations. Otherwise the conditions are not compatible.

 Let us now see that this is done in finitely many steps. The algorithm starts with
 $N_0=N$. At the next step it computes $N_1$, and then $N_2$, etc...
 These satisfy $N_0 \subseteq N_1 \subseteq N_2 \subseteq
 \cdots$. By the ACC, the process stabilizes.
 So, only a finite number of factors can exist, thus
 dropping factors is also a finite process.
 \end{proof}
The second necessary task is to reduce a given polynomial in $S$
wrt $\Sigma$. This is done in a standard form by {\tt
PNORMALFORM}. To eliminate the coefficients reducing to zero for
the given specification it suffices to compute the remainder of
the division by $N$, because $N$ is radical.
% The non-vanishing coefficients will also be reduced wrt $N$.
And then, in order to further simplify the polynomials, all those
factors lying in $W$ are also dropped from $N$. See table
\ref{pnormalformalg}.
\begin{table}
\fbox{\parbox[t]{5.2in}{

 \noindent $\tilde{f} \leftarrow \hbox{{\bf PNORMALFORM}}(f,\Sigma)$ \newline
 {\tt Input}: $f \in R[\bar{x}]$ a polynomial, $\Sigma=(N,W)$ a reduced specification, \newline
 {\tt Output}: $f$ reduced wrt $\Sigma$ \newline
 BEGIN \newline
 \hs $\widetilde{f}:=$ the product of the
 factors of $\overline{f}^N$ not laying in $W$, conveniently normalized\newline
 END

}} %end parbox, end fbox
\newline

\caption{\label{pnormalformalg}}
\end{table}

Nevertheless, the reduction using {\tt PNORMALFORM} does not
guarantee that all the coefficients of the reduced polynomial do
not cancel out for any specialization $\sigma \in \Sigma$. To test
whether adding a new coefficient to the null conditions is
compatible with $\Sigma$ we need to apply {\tt CONDTOBRANCH}.

 Given $f,g \in S$ and $\Sigma$ we say that their reduced forms $f_{\Sigma}$ and
 $g_{\Sigma}$ computed by {\tt PNORMALFORM} are equivalent wrt $\Sigma$ when
 $\sigma_{\oa}(f)$ and $\sigma_{\oa}(g)$ are proportional
 polynomials for every particular specialization $\sigma_{\oa} \in
 \Sigma$ such that $\sigma_{\oa}(\lc(f_{\Sigma})) \ne 0$ and $\sigma_{\oa}(\lc(g_{\Sigma})) \ne
 0$.

 Consider for example, $\Sigma=(N=[ab-c,ac-b,b^2-c^2],W=\phi)$,
 $f_{\Sigma}=ax+c^2$,
 $g_{\Sigma}=cx+c^2b$ and $\succ_{\oa}=\lex(a,b,c)$. $f_{\Sigma}$ and
 $g_{\Sigma}$ are not identical, but note that they are equivalent. As can be seen in this example {\tt PNORMALFORM}
 is not always able to reduce them to the same polynomial.
 Nevertheless, we have the following
 \begin{prop}
 Given two polynomials $f,g\in S$
 then $f_{\Sigma} \sim g_{\Sigma}$ wrt $\Sigma$ iff

\hs\hs\begin{minipage}[t]{5.2in}
\begin{itemize}
  \item [{\rm (i)}] $\lpp(f_{\Sigma},\succ_{\ox})=\lpp(g_{\Sigma},\succ_{\ox})$ and
  \item [{\rm (ii)}]  {\tt PNORMALFORM} applied to
  $\lc(g_{\Sigma})f_{\Sigma}-\lc(f_{\Sigma})g_{\Sigma}$ returns 0.
\end{itemize}
\end{minipage}
 \end{prop}
 \begin{proof}
 Obviously if one of both
 hypothesis fail, the reduced expressions are not equivalent wrt
 $\Sigma$.

 On the other side, suppose that (i) and (ii) hold.
 Then, using order $\succ_{\ox\oa}$ we have
 $\overline{\lc(g_{\Sigma})f_{\Sigma}}^N =
 \overline{\lc(f_{\Sigma})g_{\Sigma}}^N$ by hypothesis (ii).
 Thus, $\lc(g_{\Sigma})(\oa)\,f_{\Sigma}(\ox,\oa) =
 \lc(f_{\Sigma})(\oa)\,g_{\Sigma}(\ox,\oa)$, for all
 specializations
 in $\Sigma$. In particular it also holds for those specializations
 which do not cancel the leading coefficients of $f_{\Sigma}$ and
 $g_{\Sigma}$. And so, it follows that $f_{\Sigma}$ and $g_{\Sigma}$
 are equivalent wrt $\Sigma$.
 \end{proof}
 Thus, {\tt PNORMALFORM} does not obtain a canonical reduction of
 $f$ wrt $\Sigma$, but it can canonically recognize two equivalent reduced
 expressions.
\subsection{Reduction of Brother Final Cases with the Same $\lpp$}
In many practical computations and after applying these algorithms
to a number of cases, we have observed that some discussion trees
have pairs of terminal vertices hung from the same father vertex
with the same $\lpp$ set of their bases. As we are only interested
in those bases having different $\lpp$ sets, then each of these
brother pairs, $\{v_0,v_1\}$, can be merged in one single terminal
vertex compacting them into their father $v$ and eliminating the
distinction of the latter condition taken in $v$.

Regarding this construction, we can define a partial order
relation between two trees if, in this way, one can be transformed
into the other.
\begin{defn}
 Let $S$ and $T$ be two binary trees. We will say that $S >T$ if

\hs\hs\begin{minipage}[t]{5.2in}
 \begin{itemize}
  \item[{\rm(i)}]$T$ is a subtree of $S$ with same root and same
  intermediate vertices, and
  \item[{\rm(ii)}] for each terminal vertex $v \in T$ there is in $S$ either the
  same vertex $v \in S$ such that
  $(G_{v_T},\Sigma_{v_T})=(G_{v_S},\Sigma_{v_S})$,
  or a subtree $\overline{S} \subset S$ pending from vertex $v \in S$ with all
  its terminal vertices $u \in \overline{S}$ with $\lpp (G_ {u_ {\overline{S}}
  })= \lpp(G_{v_T})$.
  \end{itemize}
\end{minipage}
\end{defn}

 So now, given a discussion binary tree $T$, we may find the minimal tree
 $\widetilde{T}$ within the set of all trees which can be compared with
 $T$ regarding this relation. This is done by a recursive algorithm called {\tt COMPACTVERT}.

 Let us just note that the minimal tree will not have any brother terminal
 vertices with the same $\lpp$ sets of their bases.

\subsection{Rewriting the Tree with the
Discriminant Ideal} The tree $T$ built by {\tt BUILDTREE} can be
rebuilt using the discriminant ideal $N$ (see section
\ref{DiscrimIdeal}). By theorem W2, if we are given $\sigma_{\oa}$
such that there exists some $\delta \in N$ for which
$\sigma_{\oa}(\delta) \ne 0$, then $\sigma_{\oa}(I)$ corresponds
to the generic case. Thus, placing $N$ into the top vertex
labelled as $[\ ]$ in the new tree $T'$, for its non-null son
vertex we will have $T'_{[1]}=(G_{[1]},\Sigma_{[1]})$, where
$G_{[1]}$ is the generic basis and $\Sigma_{[1]}$ is a union of
specifications from $T$ corresponding to \[\Sigma_{[1]}=\{ \sigma
\ : \ \exists\, \delta\in N \hbox{ such that } \sigma(\delta) \ne
0 \} .\] No other intermediate vertices hang from this side of the
top vertex. If the strong formulation of conjecture \ref{conject}
holds, then no generic cases will hang from the first null vertex.

The subtree under the top vertex hanging from the first null son,
for which the choice is $\sigma(N)=0$, will be slightly modified
from the original $T$. The terminal vertices corresponding to
singular cases hanging from it will not be modified as, by
construction, for all of them the condition is verified by the
corresponding specifications. Thus we can rebuild the tree using
the recursive algorithm {\tt REBUILDTREE} which goes through the
old tree $T$ and rewrites the new one $T'$. At each vertex $v$ it
tests whether the condition $N$ is already included in $N_v$. If
it is the case, then it copies the whole subtree under it.
Otherwise it adds $N$ to the null ideal $N_v$ and calls {\tt
CANSPEC} to check whether the new condition is compatible or not.
If the condition is compatible then the basis will be reduced
using {\tt PNORMALFORM} and the algorithm continues. If it is not,
then the recursion stops. This algorithm produces a better new
tree with possibly less terminal cases (only generic type cases
can be dropped). This reconstruction of the tree is very little
time-consuming.

\subsection{New Generalized Gaussian Elimination {\tt GGE}}
We add here a short description of the improvements on the
generalized Gaussian elimination algorithm {\tt GGE}.

We realized, by analyzing the procedure of the old {\tt GGE}
\cite{Mo02}, that there were some special cases for which we could
guess the result of the divisions at each step and thus could be
skipped. These improvements halve the computing time.

Even though it is more efficient and faster, {\tt GGE} has become
not so useful now because the new improvements in {\tt DISPGB},
detailed above, make, in general, {\tt DISPGB} work faster without
using {\tt GGE}. So now, the use of {\tt GGE} within the execution
of {\tt DISPGB} is just optional (not used by default). However,
it can be very useful for other applications, like in the
tensegrity problem shown in section \ref{Examples}, to eliminate
some variables and simplify a given basis.

\section{Comprehensive Gr\"{o}bner Basis}\label{CGBalg}

In \cite{We02,We03} the main goal is to obtain a Comprehensive
Gr\"obner Basis. With this aim, we have built an algorithm, called
{\tt ISCGB}, to test whether a given basis $G$ is a comprehensive
Gr\"obner basis for $I$ or not. It uses {\tt PNORMALFORM}
algorithm to specialize $G$ for every terminal case in the
discussion tree. Then it checks if $\lpp(\sigma(G))$ includes the
set of $\lpp$ of the reduced Gr\"obner basis wrt $\Sigma$ for
every terminal case. If this is true for every final case then
{\tt ISCGB} returns {\tt true} otherwise returning {\tt false}.

The algorithm also informs for which cases a given basis is not a
CGB. Thus we can compute pre-images of the polynomials for which
$B$ does not specialize to a Gr\"obner basis and add them to the
given basis in order to obtain a Comprehensive Gr\"obner Basis.

 Consider a terminal case $(G_v,\Sigma_v)$  and $g\in G_v$. To simplify
 notations we do not consider the subindex $v$.
 Let $H_g=\{f_1,\dots,f_r\}$ be a basis of the ideal $I_g=I\bigcap \langle g,N
 \rangle$ whose polynomials are of the form $qg+n$, with $q\in S$ and
 $n\in \langle N \rangle$. $I_g$ contains
 all the polynomials in $I$ which can specialize to $g$ (for those
 with $\sigma(q)$ a non-null element of $R$ wrt $\Sigma$).
 Set $f'_i=\overline{f_i}^N$. Obviously, $H'_g=\{f'_1,\dots,f'_r\}$
 is a basis of $\sigma(I_g)$. Using Gr\"obner bases techniques we
 can express $g\in \sigma(I_g)$ in the form $g=\sum_i
 \alpha_if'_i$ where the $\alpha_i$'s are reduced wrt $N$, as we are
 in $I_g/N$. Then $h=\sum_i \alpha_i f_i$ specializes to $g$
 and is a pre-image of $g$ in $I$. This is used to build
 algorithm {\tt PREIMAGE} which computes a pre-image of $g$.

Combining {\tt ISCGB} and {\tt PREIMAGE} we compute a CGB using
the algorithm sketched in table \ref{cgbalg}. Let
$B=\gb(I,\succ_{\ox\oa})$, which is a tentative CGB
\cite{FoGiTr01,Ka97}, and $F=\{(G_i,\Sigma_i) \ : \ 1 \le i \le k
\}$ the set of final cases of the discussion tree built up by {\tt
DISPGB}. {\tt ISCGB} informs about the polynomials in $F$ which do
not have a pre-image in the current tentative CGB. {\tt CGB}
algorithm adds pre-images of them until a CGB is obtained.
Nevertheless, this construction is not canonical and is much more
time-consuming than building up the tree, because it uses the
product order $\succ_{\ox\oa}$ instead of working separately wrt
$\succ_{\ox}$ and $\succ_{\oa}$.
 \begin{table}
\fbox{\parbox[t]{5.2in}{
 \noindent $\tilde{B} \leftarrow \hbox{{\bf CGB}}(B,F)$ \newline
 {\tt Input}: \newline
 \hs\hs $B=\gb(I,\succ_{\ox\oa})$  \newline
 \hs\hs $F=\{(G_i,\Sigma_i) \ : \ 1 \le i \le k \}$ obtained from
 {\tt DISPGB} \newline
 {\tt Output}: $\tilde{B}$ a CGB of $I$ \newline
 BEGIN \newline
 \hs $\tilde{B}=B$ \newline
 \hs $\tilde{F}=$ SELECT cases from $F$ for which
 $\hbox{{\tt ISCGB}}(B,\succ_{\ox})$ is not a CGB. \newline
 \hs WHILE $\tilde{F}$ is non empty DO \newline
 \hs\hs TAKE the first case $(G_1,\Sigma_1)\in \tilde{F}$ \newline
 \hs\hs $\tilde{B}=\tilde{B}\ \bigcup \ \{ \hbox{{\tt PREIMAGE}}(g,\Sigma_1,B) \ : \ g \in G_1 \}$ \newline
 \hs\hs $\tilde{F}=$ SELECT cases from $\tilde{F}$ for which
 $\hbox{{\tt ISCGB}}(\tilde{B},\succ_{\ox})$ is not a CGB. \newline
 \hs END DO \newline
 END

}} %end parbox, end fbox
\newline

\caption{\label{cgbalg}}
 \end{table}

\section{Examples}\label{Examples}
% We have widely tested the implementation on
%a high number of examples from the literature.
We have selected two significative detailed examples. The first
one is the classical robot arm, which has a very nice geometrical
interpretation, and the second one is the study of a tensegrity
problem described by a linear system with the trivial null
solution in the generic case which has a non principal
discriminant ideal. After that, we outline a table containing some
relevant information for several other examples.

\subsection{Simple Robot}\label{robot}
The following system represents a simple robot arm (compare with
\cite{Mo02}):
\[
\begin{array}{l}
B=[s_1^2+c_1^2-1,s_2^2+c_2^2-1,l\,(s_1\,s_2-c_1\,c_2)-c_1+r, \\
\hs \hs l\,(s_1\,c_2+c_1\,s_2)
+s_1-z]\\
\end{array}
\]
Using the orders $\lex(s_1,c_1,s_2,c_2)$ and $\lex(r,z,l)$,
respectively for variables and parameters, {\tt DISPGB} produces
the following outputs: The discriminant ideal is principal: $N=J=
[l\,(z^2+r^2)].$ The set of final cases expressed in the form
$T_i=(G_i,(N_i,W_i))$ is:
\[
\begin{array}{lcl}
 T_{[1]} &=& ([2\,l\,c_2+l^2+1-z^2-r^2,
4\,l^2\,s_2^2+(l^2-1)^2\\
&& -2\,(l^2+1)\,(r^2+z^2)+(z^2+r^2)^2,\\
&& 2\,(r^2+z^2)\,c_1-2\,z\,l\,s_2-r\,(r^2+z^2-l^2+1), \\ &&
2\,(r^2+z^2)\,s_1+2\,l\,r\,s_2+z\,(l^2-r^2-z^2)],\ \ ([\ ], \{l\,(r^2+z^2) \} )). \\
 T_{[0, 1, 1, 1]} &=& ([2\,l\,c_2+l^2+1,
4\,(l^2-1)\,r\,c_1+2\,z\,l\,s_2-(l^2-1)\,r, \\ &&(l^2-1)^2-4\,z^2,
4\,(l^2-1)\,z\,s_1+(l^2-1)^2+4\,z^2],\\ && ([z^2+r^2], \{z, l+1,
r, l,
l-1\} )), \\
 T_{[0, 1, 1, 0]} &=& ([1],\ \ ([z, r], \{l+1, l, l-1 \})), \\
 T_{[0, 1, 0, 1]} &=& ([1],\ \ ([l^2-1,r^2+ z^2], \{z, l\} )), \\
 T_{[0, 1, 0, 0]} &=& ([l\,c_2+1, s_2, s_1^2+c_1^2-1],\ \ ([l^2-1, z, r],
\{l\} )), \\
 T_{[0, 0, 1]} &=& ([1],\ ([l],\{r^2+z^2-1\} )), \\
 T_{[0, 0, 0]} &=& ([s_2^2+c_2^2-1,c_1-r,s_1-z],\ \ ( [l,r^2+z^2-1], \{\ \})), \\
\end{array}
\]
\begin{figure}
\begin{center}
\includegraphics{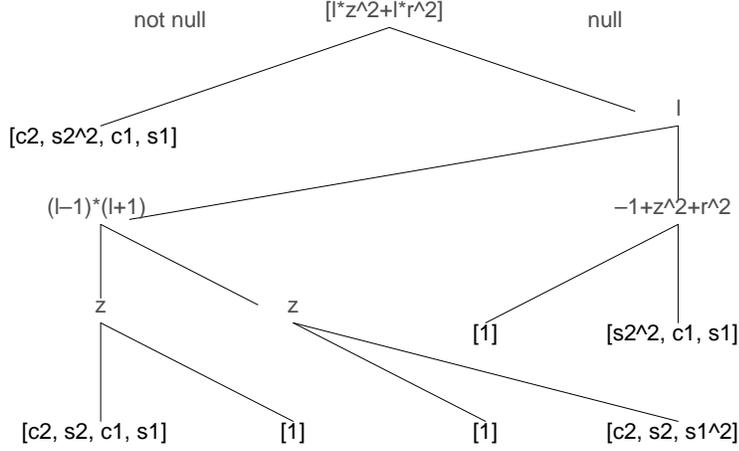}
\end{center}
\caption{\label{robotarm} {\tt DISPGB}'s graphic output for the
robot arm.}
\end{figure}

The generic case $T_{[1]}$ gives the usual formula for the robot.
%and, for this problem, is the most interesting solution.
It is characterized by the discriminant ideal $N$. The singular
cases have simple geometrical interpretation and give information
about the degenerated cases.

A graphic plot of the tree is also provided in the library. There,
the deciding conditions can be visualized at the intermediate
vertices and the $\lpp$ sets of the reduced Gr\"obner bases are
shown at the terminal vertices (see figure \ref{robotarm}).

Now we apply {\tt ISCGB} to $GB=\gb(B,\lex(s_1,c_1,s_2,c_2,r,z,l)$
wrt the output tree. The result is {\tt false}, and the list of
specializations for all the final cases is provided:

\begin{tabular}{l}
$[[1], \{s_1, s_2 c_1, s_2 s_1, c_1, c_2 s_1, c_2, s_1^2, s_2^2\},
\{s_1, c_1, c_2, s_2^2\}, \hbox{true}]]$\\
$[[0, 1, 1, 1], \{s_1, s_2, s_2 c_1, s_2 s_1, c_2 s_1, c_2, s_1^2,
s_2^2\}, \{s_1, s_2, c_1,c_2\}, \hbox{false}],$\\
$[[0, 1, 1, 0], \{1, s_1, s_2, s_2 s_1, c_2 s_1, c_2, s_1^2,
s_2^2\}, \{1\}, \hbox{true}],$ \\
$ [[0, 1, 0, 1],\{s_1, s_2, s_2 c_1, s_2 s_1, c_2 s_1, c_2, s_1^2,
s_2^2\}, \{1\},
\hbox{false}],$ \\
$[[0, 1, 0, 0], \{s_2, s_2 c_1, s_2 s_1, c_2 s_1, c_2, s_1^2,
s_2^2\}, \{s_2, c_2, s_1^2\}, \hbox{true}],$\\
$[[0, 0, 1], \{1, s_1,s_2 c_1, s_2 s_1, c_1, c_2 s_1, s_1^2,
s_2^2\}, \{1\}, \hbox{true}],$ \\
$[[0, 0, 0], \{s_1, s_2 c_1, s_2 s_1, c_1, c_2 s_1, s_1^2,
s_2^2\}, \{s_1, c_1, s_2^2\}, \hbox{true}],$ \\
\end{tabular}

There are only two cases for which $GB$ is not a CGB. Even so, the
algorithm {\tt CGB} only needs to add one single polynomial to
obtain a CGB.
\[
\begin{array}{l}
CGB = [2 l c_2+l^2+1-z^2-r^2,\  c_2^2+s_2^2-1,\  2(z^2+r^2) c_1-2
z l s_2 \\ \ \   +r(l^2-1-z^2-r^2),\  4 z s_2  c_1-4 z r s_2+4 r
c_2 c_1 +4 l r c_1\\ \ \ +2(z^2-r^2-1) c_2 -l (z^2+r^2-l^2+3),\  2
r c_1 s_2-2 z c_1 c_2
-2 z l c_1\\
\ \ +(-r^2+z^2-1+l^2) s_2+2 z r c_2,\  2(l^2-1) s_1-4 l c_1 s_2+2
l
s_2 r\\\ \  -z (r^2+z^2-l^2 -3),\  2 s_1 z+2 c_1 r-r^2-z^2+l^2-1,\\
\ \  r s_1 -z c_1+l s_2,\  s_1 c_2+ l s_1 - c_1 s_2+r s_2 -z c_2
,\
s_1 s_2+c_1 c_2\\
\ \  +l c_1-z s_2 -r c_2,\  c_1^2+s_1^2-1,\  4 (r^2+ z^2) c_1^2-4
r (1+z^2+ r^2- l^2) c_1\\ \ \  +(r^2+ z^2-l^2+1)^2-4 z^2 ].
\end{array}
\]
\subsection{Tensegrity Problem}\label{tensegrity}
We study here a problem formulated by M. de Guzm\'an and D. Orden
in \cite{GuOr04}.

Given the five points $P_1(0,0,0)$, $P_2(1,1,1)$, $P_3(0,1,0)$,
$P_4(1,0,0)$, $P_5(0,0,1)$ we want to determine a sixth one
$P_6(x,y,z)$ for which the framework with vertices
$\{P_1,\dots,P_6\}$ and edges ${{\{P_1,\dots,P_6\}} \choose {2}}
\setminus \{P_1P_6, P_2P_4, P_3P_5\}$ stays in general position
and admits a non-null self-stress.

The system describing this problem is the following:
\[ \begin{array}{l} B =
[w_{12}+w_{14}, w_{12}+w_{13}, w_{12}+w_{15},
w_{12}+w_{23}+w_{25}-w_{26} x+w_{26}, \\ \ \ w_{12}+w_{25}-w_{26}
y+w_{26}, w_{12}+w_{23}-w_{26} z+w_{26}, w_{23}+w_{34}+x w_{36},\\
\ \ w_{13}+w_{34}-w_{36} y+w_{36}, w_{23}+z w_{36},
w_{14}+w_{34}+w_{45}-w_{46} x+w_{46},\\ \ \
 w_{34}+y w_{46}, w_{45}+z
w_{56}, w_{15}+w_{45}-z w_{56}+w_{56},\\ \ \
 -w_{26}+w_{26} x+x
w_{36}-w_{46}+w_{46} x+w_{56} x, \\ \ \ -w_{26}+w_{26}
y-w_{36}+w_{36} y+y w_{46}+w_{56} y,
\\ \ \ -w_{26}+w_{26} z+z w_{36}+w_{46} z-w_{56}+z w_{56}]
 \end{array} \]
\begin{figure}
\begin{center}
\includegraphics{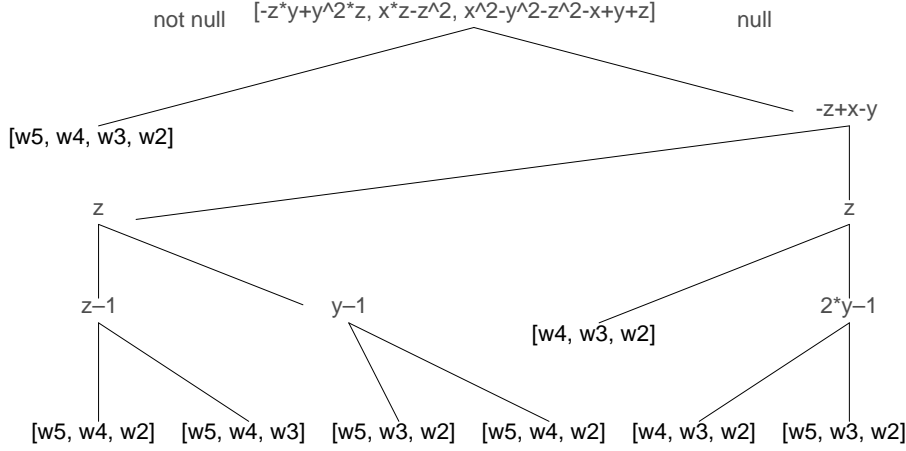}
\end{center}
\caption{\label{tensegritydraw} {\tt DISPGB} graphic output for
the tensegrity problem.}
\end{figure}
Set
$\succ_{\ox\oa}=\lex(w_{12},w_{13},w_{14},w_{15},w_{23},w_{25},w_{34},w_{45},
w_{26},w_{36},w_{46},w_{56},x,y,z)$. In order to simplify the
system we compute $\hbox{{\tt GGE}}(B,\succ_{\ox\oa})$
(Generalized Gaussian Elimination). The {\tt GGE} basis can be
expressed as $B'=B'_1 \cup B'_2$, with $B'_2$ being the
elimination ideal wrt the variables $w_{26}=w_{2}$, $
w_{36}=w_{3}$, $w_{46}=w_{4}$, $w_{56}=w_{5}$, and $B'_1$
expressing the remaining variables linearly in terms of
$w_{2},w_{3},w_{4},w_{5}$:
\[ \begin{array}{l}
B'_1=[w_{45}+z w_{5}, w_{34}+y w_{4}, w_{25}+w_{5} y, w_{23}+z
w_{3}, w_{15}-2 z w_{5}+w_{5},\\ \ \ \   w_{14}-2 z w_{5}+w_{5},
w_{13}-2 z w_{5}+w_{5}, w_{12}+2 z w_{5}-w_{5}] \\
B'_2 = [-z w_5+w_5 x-w_5 y, -z w_5+w_4 z, w_4 x+y w_4-w_4-z
w_5+w_5,\\ \ \ \  w_3 y-w_3+y w_4-2 z w_5+w_5, x w_3-y w_4-z
w_3,\\ \ \ \  w_2 z-w_2+z w_3+2 z w_5-w_5, w_2 y-w_2+w_5 y+2 z
w_5-w_5, \\ \ \ \ w_2 x-w_2+z w_3+w_5 y+2 z w_5-w_5],\\
 \end{array} \]
Then, using the orders $\succ_{\ox}=\lex(w_2,w_3,w_4,w_5)$ and
$\succ_{\oa}=\lex(x,y,z)$, for variables and parameters
respectively, {\tt DISPGB} produces the following output:

\[\begin{array}{lcl}
 T_{[1]} &=& ([w_5, w_4, w_3, w_2],
([\ ], \\ & & \{[y^2 z-y z, z x-z^2, x^2-y^2-z^2-x+y+z]\}))\\
 T_{[0, 1, 1, 1]} &=& ([w_5, w_4, w_2 z-w_2+z w_3],
([y-1, x-z],\{z, z-1\}), \\
 T_{[0, 1, 1, 0]} &=& ([w_5, w_4, w_3], ([z-1, y-1, x-1], \{\
\})), \\
 T_{[0, 1, 0, 1]} &=& ([w_5,
y w_4+w_3 y-w_3, w_2], ([z, y-1+x], \{2 y-1, y-1\}),\\
 T_{[0, 1, 0, 0]} &=& ([w_5, w_4, w_2], ([z, y-1, x], \{\ \})),\\
 T_{[0, 0, 1]} &=& ([-w_5+w_4,
w_3+2 z w_5-w_5, w_2-2 z w_5+w_5],\\ && ([y, x-z],\{z\})),\\
\end{array}
\]
\[\begin{array}{lcl}
 T_{[0, 0, 0, 1]} &=& ([w_5+2 y w_4-w_4, 2 w_3 y-w_3+w_5, w_2+w_5],
 \\ && ([z, x-y], \{2 y-1\})), \\
 T_{[0, 0, 0, 0]} &=& ([w_5, w_3-w_4, w_2], ([z, 2 y-1, 2 x-1], \{\
\})),\\
\end{array}\]
\begin{figure}
\begin{center}
\includegraphics{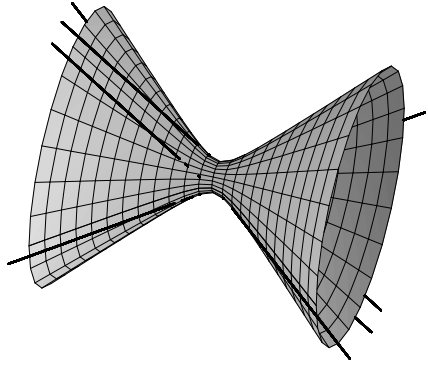}
\end{center}
\caption{\label{tensegritydrawA} Location of the sixth point for
non null self-stress. }
\end{figure}
and the discriminant ideal is not principal:
\[N=J= [y^2 z-y z, z x-z^2, x^2-y^2-z^2-x+y+z].\]

The generic solution is trivial ($w_5= w_4= w_3= w_2= 0$). In this
problem, the interesting non trivial solutions are given by the
conditions over the parameters described by the variety of the
discriminant ideal, which decomposes into 4 straight lines
included in the hyperboloid $x^2-y^2-z^2-x+y+z=0$ (illustrated in
figure \ref{tensegritydrawA}):
\[ \V(N)=\V(z,x-y)\bigcup \V(y,x-z) \bigcup \V(z,x+y-1) \bigcup
\V(y-1,x-z).\]

For this problem the Gr\"obner basis wrt variables and parameters
is already a comprehensive Gr\"obner basis.

\section{Benchmarks} For a set of examples taken from the
literature we have applied the current implementation, release 2.3
in {\em Maple 8 }, of algorithm {\tt DISPGB} using a 2 GHz Pentium
4 at 512 MB. Table \ref{bench} summarizes the computing time of
{\tt DISPGB}, the total number of terminal vertices of the output
tree, whether the discriminant ideal is principal or not, and
whether the D-Gr\"{o}bner basis wrt $\succ_{\ox\oa}$ is a CGB or not,
joint by the number of failure cases for which it is not (0 if it
is). The bases of the different examples are detailed below:
\begin{itemize}
\item {S1.} $[a (x+y), b (x+y), x^2+a x]$;
\item {S2.}
$[x_1^2,x_1 x_2,x_1 x_3^2,x_1 a+x_2,x_2 x_3-x_3^2,x_2
a,x_3^3,x_3^2 a,a^2]$;
\item{S3.} $[x^3-a x y,x^2 y-2 y^2+x]$;
\item {S4.} $[a x+y-1, b x+y-2, 2 x+a y, b x+a y+1]$;
\item {S5.}
$[x_4-(a_4-a_2),x_1+x_2+x_3+x_4-(a_1+a_3+a_4),  x_1 x_3 x_4-a_1
a_3 a_4,\\ \null \hspace{8mm} x_1 x_3+x_1 x_4+x_2 x_3+x_3 x_4-(a_1
a_4+a_1 a_3+a_3 a_4)]$;
\item {S6.} $[v x y+u x^2+x,u y^2+x^2]$;
\item {S7.}
 $[y^2-z x y+x^2+z-1,x y+z^2-1,y^2+x^2+z^2-r^2]$;
\item{S8.}
$[a-b+(x y a-x^2 y b-3 a)^3+(x y b-3 x b-5 b)^4,x y a-x^2 y b-3
a,\\ \null \hspace{9mm} x y b-3 x b-5 b]$;
\item {S9.} $[x+c y+b z+a, c x+y+a z+b, b x+a y+z+c]$;
\item {S10. See subsection \ref{robot}};
\item{S11.} $[(d_4 d_3 R+r_2^2-d_4 d_3 r_2^2+d_4^2 d_3^2-d_4 d_3^3-d_4^3 d_3+d_4 d_3+Z-R) t^4\\
\null \hspace{9mm} +(-2 r_2 d_4 R+2 r_2 d_4^3+2 r_2 d_4 d_3^2-4 r_2 d_3 d_4^2+2 r_2^3 d_4+2 r_2 d_4) t^3\\
\null \hspace{9mm}
-(2 r_2^2-2 R+4 d_4^2 r_2^2+4 d_4^2+2 Z-2 d_4^2 d_3^2) t^2\\
\null \hspace{9mm}
+(-2 r_2 d_4 R+2 r_2 d_4 d_3^2+2 r_2 d_4+2 r_2 d_4^3+4 r_2 d_3 d_4^2+2 r_2^3 d_4) t\\
\null \hspace{9mm} +r_2^2+d_4^3 d_3-d_4 d_3 R+d_4 d_3
r_2^2+Z-R-d_4 d_3+d_4^2 d_3^2+d_4 d_3^3]$;
\item{S12.}
$[a-l_3 c_3-l_2 c_1,b-l_3 s_3-l_2
s_1,c_1^2+s_1^2-1,c_3^2+s_3^2-1]$;
\item{S13.}
$[a x^2 y+a+3 b^2,a (b-c) x y+a b x+5 c]$;
\item {S14.} $[t^3-c u t^2-u v^2-u w^2, t^3-c v t^2-v u^2-v w^2,
t^3-c w t^2-w u^2-w v^2]$;
\item {S15.} $ [a+d s_1,b-d c_1,l_2 c_2+l_3 c_3-d, l_2 s_2+l_3 s_3-c_,
 s_1^2+c_1^2-1,s_2^2+c_2^2-1,\\ \null \hspace{9mm}
 s_3^2+c_3^2-1]$;
\item {S16. See subsection \ref{tensegrity}.}
\end{itemize}
We have tested several other problems and in some of them only
partial results have been reached. We detail two significative
examples:
\begin{itemize}
\item{S17.} $[a x t^2+b y t z-x (x^2+c y^2+d z^2), a y t^2+b z x t-x (y^2+c
z^2+d x^2), \\  \null \hspace{9mm}
a z t^2+b x y t-x (z^2+c x^2+d
y^2)]$
\item{S18.} $
[(3 x^2+9 v^2-3 v-3 x) t_1^2 t_2^2 +(3 v^2-3 v+6 v x-3 x+3 x^2) t_2^2 \\
 \null \hspace{9mm}+(3 v+3 v^2+3 x^2-3 x-6 v x) t_1^2-24 v^2 t_1 t_2+9 v^2-3 x+3 x^2+3 v, \\
 \null \hspace{9mm}(3 x^2+9 v^2-3 v-3 x) t_2^2 t_3^2+(3 v+3 v^2+3 x^2-3 x-6 v x) t_2^2 \\
 \null \hspace{9mm}+(3 v^2-3 v+6 v x-3 x+3 x^2) t_3^2-24 v^2 t_2 t_3+9 v^2-3 x+3 x^2+3 v,\\
 \null \hspace{9mm} (3 x^2+9 v^2-3 v-3 x) t_3^2 t_1^2 +(3 v^2-3 v+6 v x-3 x+3 x^2) t_1^2\\
 \null \hspace{9mm}+(3 v+3 v^2+3 x^2-3 x-6 v x) t_3^2-24 v^2 t_3 t_1+9 v^2-3 x+3
x^2+3 v]$
\end{itemize}
For S17 \cite{GoTrZa00}, {\tt DISPGB} gets bogged down after
computing 35 terminal vertices in 1375 sec. It has been unable to
finish the tree, and so neither rebuilding with the
discriminantideal nor reducing the tree can have been achieved.
The label of the 35th vertex is $[1,1,0,1,0,0]$, thus all vertices
beginning with $[0,0,\dots$ have been already determined (the tree
is built up in pre-order beginning with the 0 vertices).

S18 corresponds to the benzene molecule studied in \cite{Em99}.
The situation is similar to S17, getting bogged down after 45
seconds when the 9th vertex labelled  $[1,1,0,0]]$ has been
computed.
\begin{table}[t]
\begin{center}
\begin{tabular}{||l|c|c|c|c||}
\hline\hline Identification & CPU time      & \# Final & Discriminant  & Is CGB? \\
                  & (seconds) & vertices & is principal? & (\# failures) \\
                  \hline\hline
\hline S1 \cite{We03} & 0.8 & 2 & N & Y (0) \\
\hline S2. \cite{Gi87} & 1.2 & 2 & Y & N (1) \\
\hline S3. \cite{Gom02,Du95} & 1.5 & 2 & Y & Y (0) \\
\hline S4. & 1.6 & 2 & N & Y (0) \\
\hline S5. \cite{Kap95} & 1.6 & 3 & Y & N (1) \\
\hline S6. \cite{Kap95} & 2.0 & 4 & Y & Y (0) \\
\hline S7. \cite{Kap95} & 3.0 & 2 & Y & Y (0) \\
\hline S8. \cite{SaSuNa03} & 4.4 & 3 & Y & Y (0) \\
\hline S9. Similar to \cite{Si92} & 6.7 & 10 & Y & Y (0) \\
\hline S10. Subsection \ref{robot} & 7.9 & 7 & Y & N (2) \\
\ \ \ \ \ \ Simple robot & & & & \\
\hline S11. \cite{Co04} & 8.0 & 6 & Y & Y (0) \\
\ \ \ \ \ \ Singular points & & & & \\
\hline S12. \cite{Ry00} & 8.2 & 11 & Y & N (1) \\
\ \ \ \ \ \ Rychlik robot & & & & \\
\hline S13. \cite{SaSu03} & 8.2 & 10 & Y & Y (0) \\
\hline S14. \cite{GoTrZa00,De99} & 9.6 & 2 & Y & N (1) \\
\ \ \ \ \ \ ROMIN robot & & & & \\
\hline S15. \cite{GoRe93} & 18.2 & 17 & Y & N (2) \\
\hline S16. \cite{GuOr04} & 21.3 & 8 & N & Y (0) \\
\ \ \ \ \ \ Subsection \ref{tensegrity} & & & & \\
 \hline\hline
\end{tabular}
\newline

\caption{\label{bench}}
\end{center}
\end{table}
\section{Acknowledgements}
We want to thank Professor Volker Weispfenning for his useful
suggestions and ideas and for encouraging us in undertaking the
research on Comprehensive Gr\"obner Bases, as well as for his
hospitality and kindness on the occasions in which we met.

We would like to thank Professor Pelegr\'{\i} Viader for his many
helpful comments and his insightful perusal of our first draft.

We will also thank the referees for their valuable suggestions.

\bibliographystyle{elsart-harv}

\end{document}